\documentclass[12pt,letterpaper,titlepage]{amsart}
\usepackage{amsmath, amssymb, amsthm, amsfonts,amscd,xr}
\newcommand{\N}{\mathbb{N}}

\newcommand{\Z}{\mathbb{Z}}
\newcommand{\R}{\mathbb{R}}

\def\beq{\begin{equation}}
\def\eeq{\end{equation}}

\def\arr{\hbox to 20pt{\rightarrowfill}}

\def\d{\delta}
\def\e{\varepsilon}

\def\Sl{\mathrm {Sl}}

\def\U{\mathcal U}

\newenvironment{res} 
               {\begin{equation} 
\begin{minipage}{0.85\textwidth}} 
               { \end{minipage}\end{equation} } 
\def\ber{\begin{res} } 
\def\eer{\end{res}} 
 
\numberwithin{equation}{section} 
\newtheorem{thm}{Theorem}[section]

\makeatletter 
\def\section{\@startsection {section}{1}{\z@}{3.5ex plus 1ex minus 
    .2ex}{2.3ex plus .2ex}{\large\bf}} 
    \def\subsection{\@startsection{subsection}{2}{\z@}{3.25ex plus 1ex minus 
 .2ex}{1.5ex plus .2ex}{\bf}} 
\makeatother 
\def\pf{{\em Proof}.\, }

\def\ad{\operatorname{ad}}

\def\e{\epsilon}

\def\d{\delta}

\def\af{\mathfrak{a}} 

\def\gf{\mathfrak{g}} 
\def\h{\mathfrak{h}}

\def\kf{\mathfrak{k}}

\def\pf{\mathfrak{p}} 
 
\def\qf{\mathfrak{q}}

\def\bs{\backslash}

\makeindex 
\begin{document} 
\title[Function spaces]
{On function spaces on symmetric spaces}
\author{Bernhard Kr\"otz} 
\address{Max-Planck-Institut f\"ur Mathematik\\  
Vivatsgasse 7\\ D-53111 Bonn
\\email: kroetz@mpim-bonn.mpg.de}

\author{Henrik Schlichtkrull}
\address{Department of Mathematics\\ 
University of Copenhagen\\ Universitetsparken 5 \\ 
DK-2100 K\o{}benhavn 
\\ email: schlicht@math.ku.dk} 

\date{\today} 
\begin{abstract}{Let $Y=G/H$ be a semisimple symmetric space.
It is shown that the smooth vectors for the regular representation
of $G$ on $L^p(Y)$ vanish at infinity.}\end{abstract}
\thanks{}
\maketitle 

\section{Vanishing at infinity}

Let $G$ be a connected unimodular Lie group, equipped with a
Haar measure $dg$, and let $1\leq p<\infty$. 
We consider the left regular representation $L$
of $G$ on the function space $E_p=L^p(G)$.

Recall that $f\in E_p$ is called a 
{\it smooth vector for} $L$ if and only if 
the map 
$$G\to E_p, \ \ g\mapsto L(g) f$$
is a smooth $E_p$-valued map. 

Write $\gf$ for the Lie algebra of $G$ and 
$\U(\gf)$ for its enveloping algebra. 
The following result is well-known, see \cite{Poulsen}.

\begin{thm}\label{result for group} 
The space of smooth vectors for $L$
is
$$E_p^\infty=\{f\in C^\infty(G)\mid L_uf\in L^p(G)
\text{ for all } u\in\U(\gf)\}.
$$
Furthermore, $E_p^\infty\subset C_0^\infty(G)$,
the space of smooth functions
on $G$ which vanish at infinity.
\end{thm}

Our concern is with the corresponding result
for a homogeneous space $Y$ of $G$. By that we 
mean a connected manifold $Y$ with a transitive action of 
$G$. In other words 
$$Y=G/H$$
with $H\subset G$ a closed subgroup. We shall request 
that $Y$ carries a $G$-invariant positive measure 
$dy$. Such a measure is unique up to scale and 
commonly referred to as Haar measure. 
With respect to $dy$ we form the Banach spaces 
$E_p:=L^p(Y)$. The group 
$G$ acts continuously by isometries on $E_p$ via 
the left regular representation:
$$[L(g)f](y)=f(g^{-1}y)\qquad (g\in G, y\in Y, f\in E_p)\, .$$
We are concerned with the space $E_p^\infty$
of smooth vectors 
for this representation. The first part of
Theorem \ref{result for group}
is generalized as follows, see 
\cite{Poulsen}, Thm. 5.1.

\begin{thm}
The space of smooth vectors for $L$
is
$$E_p^\infty=\{f\in C^\infty(Y)\mid L_uf\in L^p(Y)
\text{ for all } u\in\U(\gf)\}.
$$
\end{thm}

We write $C_0^\infty(Y)$ for the space of smooth functions
vanishing at infinity. Our goal 
is to investigate an assumption  under which
the second part of Theorem \ref{result for group}
generalizes, that is,
\begin{equation} 
\label{v} E_p^\infty\subset C_0^\infty(Y).
\end{equation} 

Notice that if $H$ is compact, then we can regard
$L^p(G/H)$ as a closed $G$-invariant subspace of $L^p(G)$,
and (\ref{v}) follows immediately from Theorem
\ref{result for group}.

Likewise, if $Y=G$ regarded as a homogeneous space
for $G\times G$ with the left$\times$right action,
then again (\ref{v}) follows from Theorem
\ref{result for group}, since a left$\times$right
smooth vector is obviously also left smooth.

However, (\ref{v}) is false in general as the following class 
of examples shows. Assume that $Y$ has finite volume but is 
not compact, e.g. $Y=\Sl(2,\R)/ \Sl(2,\Z)$. Then the constant 
function ${\bf 1}_Y$ is a smooth vector for $E^p$,
but it does not vanish at infinity. 

\section{Proof by convolution}
We give a short proof of (\ref{v}) for the case $Y=G$,
based on the theorem of Dixmier and Malliavin 
(see \cite{DM}). According to this theorem,
every smooth vector in
a Fr\'echet representation $(\pi,E)$ belongs to the G\aa rding space,
that is, it is spanned by vectors of the form
$\pi(f)v$, where $f\in C^\infty_c(G)$ and $v\in E$.
Let such a vector $L(f)g$, where $g\in E_p=L^p(G)$, be given.
Then by unimodularity
\begin{equation}\label{convolution identity}
[L(f)g](y)=\int_G f(x)g(x^{-1}y)\,dx =
\int_G f(yx^{-1})g(x)\,dx.
\end{equation}
For simplicity we assume $p=1$. The general case is similar.
Let $\Omega\subset G$ be compact such that
$|g|$ integrates to $<\epsilon$ over the complement. 
Then, for $y$ outside of
the compact set ${\rm supp}f\cdot\Omega$, we have
$$yx^{-1}\in{\rm supp}f\Rightarrow x\notin\Omega,$$
and hence
$$|L(f)g(y)|\leq \sup|f|\int_ {x\notin\Omega} |g(x)|\,dx
\leq\sup|f|\,\epsilon.$$
It follows that $L(f)g\in C_0(G)$.

Notice that the assumption $Y=G$ is crucial in this proof, since the
convolution identity (\ref{convolution identity})
makes no sense in the general case.

\section{Semisimple symmetric spaces}

Let $Y=G/H$ 
be a semisimple symmetric space.
By this we mean: 
\begin{itemize}
\item $G$ is a connected semisimple Lie group with 
finite center. 
\item There exists an involutive automorphism $\tau$ of $G$ 
such that $H$ is an open subgroup of the group 
$G^\tau=\{ g\in G\mid \tau(g)=g\}$ of 
$\tau$-fixed points. 
\end{itemize}
We will verify (\ref{v}) for this case. In fact, our proof
is valid also under the more general assumption that $G/H$
is a reductive symmetric space of Harish-Chandra's class,
see \cite{Ban}.

\begin{thm} Let $Y=G/H$ be a semisimple symmetric space, 
and let $E_p=L^p(Y)$ where $1\leq p<\infty$. Then 
$$E_p^\infty\subset C_0^\infty(Y)\, .$$
\end{thm}

\begin{proof} 
A little bit of standard terminology is useful. As customary
we use the same symbol for an automorphism of $G$ and its 
derived automorphism of the Lie algebra $\gf$. 
Let us write $\gf=\h +\qf$ for the decomposition 
in $\tau$-eigenspaces according to eigenvalues 
$+1$ and $-1$.  

\par Denote by $K$ a maximal compact subgroup of $G$. We will 
and may assume that $K$ is stable under $\tau$. Write $\theta$ 
for the Cartan-involution on $G$ with fixed point group 
$K$ and write $\gf=\kf+\pf$ for the eigenspace decomposition 
of $\theta$. We fix a maximal abelian 
subspace $\af\subset\pf\cap \qf$. 

\par The simultaneous eigenspace decomposition 
of $\gf$ under $\ad \af$ leads to a (possibly reduced) root system
$\Sigma\subset \af^*\bs\{0\}$. 
Write $\af_{\rm reg}$ for $\af$ with the root hyperplanes 
removed, i.e.: 
$$\af_{\rm reg}=\{ X\in\af\mid (\forall \alpha\in\Sigma)\ \alpha(X)\neq 0\}\, .$$
Let $M=Z_{H\cap K}(\af)$ 
and $W_H=N_{H\cap K}(\af)/ M $. 
%for the small Weyl group of $\Sigma$

Recall the polar decomposition of $Y$.
With $y_0=H\in Y$ the base point of $Y$ it asserts that the mapping 
\def\polar{\rho}
$$\polar: K/M \times \af \to Y, \ \ (kM, X)\mapsto k\exp(X)\cdot y_0$$
is differentiable, onto and proper. Furthermore,
the element $X$ in the decomposition is unique up
to conjugation by $W_H$, and the induced map 
$$K/M\times_{W_H}\af_{\rm reg} \to Y$$
is a diffeomorphism onto an open and dense subset of $Y$.

Let us return now to our subject proper, the vanishing 
at infinity of functions in $E_p^\infty$. 
Let us denote functions on $Y$ by lower case roman letters,
and by the corresponding upper case letters
their pull backs to $K/M\times \af$, for example $F=f\circ\polar$. 
Then $f$ vanishes at infinity on $Y$ translates into 
\begin{equation}\label{vi} \lim_{{X\to\infty}\atop{X\in \af}} 
\sup_{k\in K} |F(kM, X)|=0\, .\end{equation}

\par We recall the formula 
for the pull back by $\polar$ of the invariant measure $dy$ on $Y$. 
For each $\alpha\in\Sigma$ we 
denote by $\gf^\alpha\subset\gf$ the corresponding root space. 
We note that $\gf^\alpha$ is stable under the involution $\theta\tau$. 
Define $p_\alpha$, resp. $q_\alpha$, as the dimension of the 
$\theta\tau$-eigenspace in $\gf^\alpha$ according to 
eigenvalues $+1,-1$. 
Define a function $J$ on $\af$ by 
$$J(X)=\left| \prod_{\alpha\in\Sigma^+} [\cosh \alpha(X)]^{q_\alpha}\cdot 
[\sinh\alpha(X)]^{p_\alpha}\right|\, .$$  

With $d(kM)$ the Haar-measure on $K/M$ and $dX$ the Lebesgue-measure 
on $\af$ one then gets, up to normalization: 
$$\polar^*(dy) = J(X)\, d(k,X):=J(X)\  d(kM)\, dX\, .$$ 

We shall use this formula to
relate certain Sobolev norms on $Y$ and on $K/M\times \af$.
Fix a basis 
$X_1, \ldots, X_n$ for $\gf$. For an $n$-tupel 
${\bf m}=(m_1, \ldots, m_n)\in \N_0^n$ we define 
elements $X^{\bf m}\in \U(\gf)$ by 
$$X^{\bf m}:=X_1 ^{m_1}\cdot \ldots\cdot X_n^{m_n}\, .$$   
These elements form a basis for $\U(\gf)$. 
We introduce the $L^p$-Sobolev norms on $Y$,
$$S_{m,\Omega}(f):=
\sum_{|{\bf m}|\leq m} \left[\int_\Omega
|L(X^{\bf m})f(y)|^p\,dy\right]^{1/p}$$
where $\Omega\subset Y$, and where
$|{\bf m}|:=m_1+\ldots+ m_n$. Then $f\in E_p^\infty$
if and only if $S_{m,Y}(f)<\infty$ for all $m$.

Likewise, for $V\subset\af$ 
we denote 
$$S^*_{m,V}(F):=
\sum_{|{\bf m}|\leq m} \left[\int_{K\times V}
|L(Z^{\bf m})F(kM,X)|^p \,
J(X)\,d(k,X)\right]^{1/p}$$
Here $Z$ refers to members of some fixed
bases for $\kf$ and $\af$, acting from the left
on the two variables,
and again $\bf m$ is a multiindex.

Observe that for $Z\in\af$ we have for the action
on $\af$,
$$[L(Z)F](kM,X)=[L(Z^k)f](k\exp(X)\cdot y_0)$$
where $Z^k:={\rm Ad}(k)(Z)$ can be written 
as a linear combination of the basis elements in $\gf$, 
with coefficients which are continuous on $K$.
It follows that there exists a constant
$C_m>0$ such that for all $F=f\circ\polar$,
\begin{equation}
\label{sk}
S^*_{m,V}(F)\leq C_{m} S_{m,\Omega}(f)
\end{equation}
where $\Omega=\polar(K/M,V)=K\exp(V)\cdot y_0$.

Let $\e>0$ and set 
$$\af_\e:=\{X\in\af\mid (\forall \alpha\in\Sigma)\ 
|\alpha(X)|\geq \e\}\, . $$
Observe that there exists a constant $C_\e>0$ such that 
\begin{equation}\label{eps}(\forall X\in \af_\e) 
\quad J(X)\geq C_\e \,. \end{equation}
 
We come to the main part of the proof.
Let $f\in E_p^\infty$.
We shall first establish that
\begin{equation}\label{fs}  
\lim_{X\to \infty\atop X\in \af_\e} 
F(eM, X)=0\,.
\end{equation} 

It follows from the Sobolev lemma,
applied in local coordinates,
that the following holds for a sufficiently large 
integer $m$ (depending only on $p$ and the dimensions 
of $K/M$ and $\af$).
For each compact symmetric neighborhood $V$ of $0$ in
$\af$ there exists a constant $C>0$ such that
\begin{equation}\label{Sobolev}
\begin{aligned}
&|F(eM,0)| \\&\quad\leq
C \sum_{|{\bf m}|\leq m}
\left[
\int_{K/M\times V} |[L(Z^{\bf m})F](kM,X)|^p \,
d(k,X)
\right]^{1/p}
\end{aligned}
\end{equation}
for all $F\in C^\infty(K/M\times\af)$.
We choose $V$ such that
$\af_\e+V\subset\af_{\e/2}$.

Let $\d>0$. Since  $f\in E^p$, 
it follows from (\ref{sk}) and the properness of $\polar$ that
there exists a compact set $B\subset\af$ with complement
$B^c\subset\af$,
such that
\begin{equation}
\label{defi B}
S^*_{m,B^c}(F)\leq C_m S_{m,\Omega}(f)<\delta
\end{equation}
where $\Omega=K\exp(B^c)\cdot y_0$.

Let  $X_1\in \af_\e\cap(B+V)^c$. Then
$X_1+X\in\af_{\e/2}\cap B^c$ for $X\in V$.
Applying (\ref{Sobolev}) to the function
$$F_1(kM,X)=F(kM,X_1+X),$$
and employing (\ref{eps}) for the set $\af_{\e/2}$,
we derive
$$\begin{aligned}
&\,|F(eM,X_1)|\\
&\quad\leq
C\sum_{|{\bf m}|\leq m}
\left[
\int_{K/M\times V}
 |[L(Z^{\bf m})F_1](kM,X)|^p \,
d(k,X)
\right]^{1/p}\\
&\quad\leq C'
\sum_{|{\bf m}|\leq m}
\left[
\int_{K/M\times B^c}
 |[L(Z^{\bf m})F](kM,X)|^p \,
J(X)\,d(k,X)
\right]^{1/p}\\
&\quad=C'S^*_{m,B^c}(F)\leq C'\d,
\end{aligned}
$$
from which (\ref{fs}) follows.

In order to conclude the theorem,
we need a version of (\ref{fs}) which is uniform for all 
functions $L(q)f$, for $q\in Q\subset G$ a compact subset. 

Let $\delta>0$ be given, and as before let $B\subset\af$ be
such that (\ref{defi B}) holds. 
By the properness of $\polar$,
there exists a compact set $B'\subset\af$ such that
$$QK\exp(B)\cdot y_0\subset K\exp(B')\cdot y_0.$$

We may assume that $B'$ is $W_H$-invariant. Then, for
each $k\in K$, $X\notin B'$ and $q\in Q$ we have that
\begin{equation}
\label{Q-orbit}
q^{-1}k\exp(X)\cdot y_0\notin K\exp(B)\cdot y_0,
\end{equation}
since otherwise we would have
$$k\exp(X)\cdot y_0\in qK\exp(B)\cdot y_0
\subset K\exp(B')\cdot y_0$$
and hence $X\in B'$. 

We proceed as before, with $B$ replaced by $B'$,
and with $f$, $F$ replaced
by $f_q=L_{q}f$, $F_q=f_q\circ p$.
We thus obtain for $X_1\in \af_\e\cap(B'+V)^c$,
$$
|F_q(eM,X_1)|\leq C S^*_{m,(B')^c}(F_q)
\leq C C_m S_{m,\Omega'}(f_q)
$$
where $\Omega'=K\exp((B')^c)\cdot y_0$.

Observe that for each $X$ in $\gf$ 
the derivative $L(X)f_q$ can be written as a linear combination 
of derivatives of $f$ by basis elements from
$\gf$, with coefficients which are uniformly bounded on $Q$.
We conclude that $S_{m,\Omega'}(f_q)$
is bounded by a constant times
$S_{m,Q^{-1}\Omega'}(f)$, with a uniform
constant for $q\in Q$. By (\ref{Q-orbit}) and
(\ref{defi B}) we 
conclude that the latter Sobolev norm is bounded
from the above by $\delta$.
 
We derive the desired uniformity of the limit 
(\ref{fs}) for $q\in Q$,
\begin{equation}\label{fs1}  
\lim_{X\to \infty\atop X\in \af_\e}\sup_{q\in Q}\, 
|F_q(eM, X)|=0\, .\end{equation}

\par Finally we choose an appropriate
set $Q$. Let $\epsilon>0$
be arbitrary.  
There exists $X_1, \ldots, X_N \in \af$ such that 
\begin{equation} \label{N} \af=\bigcup_{j=1}^N (X_j+ \af_\e )\, .
\end{equation} 
Set $a_j=\exp(X_j)\in A$ and define a compact subset of $G$ by 
$$Q:=\bigcup_{j=1}^N Ka_j \, .$$ 
Then, for every $X\in\af$ we have
$X-X_j\in\af_\e$ for some $j$. Hence with $q=k\exp(X_j)$
$$\lim_{X\to \infty} F(kM,X)=
\lim_{X\to \infty} F_q(eM,X-X_j)= 0,$$
as was to be shown. 
\end{proof}

\noindent{\bf Remark.}  Let $f\in L^2(Y)$ be a $K$-finite 
function which is also finite for the center of $\U(\gf)$.
Then it follows from \cite{RS} that $f$ vanishes at infinity. 
The present result is more general, since such a function
necessarily belongs to $E_2^\infty$.

\end{document}